\documentclass[reqno, 12pt]{amsart}
\pdfoutput=1
\makeatletter
\let\origsection=\section \def\section{\@ifstar{\origsection*}{\mysection}}
\def\mysection{\@startsection{section}{1}\z@{.7\linespacing\@plus\linespacing}{.5\linespacing}{\normalfont\scshape\centering\S}}
\makeatother

\usepackage{amsmath,amssymb,amsthm}
\usepackage{mathrsfs}
\usepackage{mathabx}\changenotsign
\usepackage{dsfont}
\usepackage{bbm}

\usepackage{xcolor}
\usepackage[backref=section]{hyperref}
\usepackage[ocgcolorlinks]{ocgx2}
\hypersetup{
	colorlinks=true,
	linkcolor={red!60!black},
	citecolor={green!60!black},
	urlcolor={blue!60!black},
}

\usepackage[open,openlevel=2,atend]{bookmark}

\usepackage[abbrev,msc-links,backrefs]{amsrefs}
\usepackage{doi}

\renewcommand{\PrintDOI}[1]{\doi{#1}}

\usepackage[T1]{fontenc}
\usepackage{lmodern}
\usepackage[babel]{microtype}
\usepackage[english]{babel}

\usepackage{geometry}
\numberwithin{equation}{section}
\numberwithin{figure}{section}

\usepackage{enumitem}

\let\polishlcross=\l
\def\l{\ifmmode\ell\else\polishlcross\fi}

\let\emptyset=\varnothing
\let\setminus=\smallsetminus

\makeatletter
\def\moverlay{\mathpalette\mov@rlay}
\def\mov@rlay#1#2{\leavevmode\vtop{   \baselineskip\z@skip \lineskiplimit-\maxdimen
		\ialign{\hfil$\m@th#1##$\hfil\cr#2\crcr}}}
\newcommand{\charfusion}[3][\mathord]{
	#1{\ifx#1\mathop\vphantom{#2}\fi
		\mathpalette\mov@rlay{#2\cr#3}
	}
	\ifx#1\mathop\expandafter\displaylimits\fi}
\makeatother

\newcommand{\dcup}{\charfusion[\mathbin]{\cup}{\cdot}}

\DeclareFontFamily{U}  {MnSymbolC}{}
\DeclareSymbolFont{MnSyC}         {U}  {MnSymbolC}{m}{n}
\DeclareFontShape{U}{MnSymbolC}{m}{n}{
	<-6>  MnSymbolC5
	<6-7>  MnSymbolC6
	<7-8>  MnSymbolC7
	<8-9>  MnSymbolC8
	<9-10> MnSymbolC9
	<10-12> MnSymbolC10
	<12->   MnSymbolC12}{}
\DeclareMathSymbol{\powerset}{\mathord}{MnSyC}{180}

\usepackage{tikz}
\usetikzlibrary{calc,decorations.pathmorphing}
\usetikzlibrary{arrows,decorations.pathreplacing}
\pgfdeclarelayer{background}
\pgfdeclarelayer{foreground}
\pgfdeclarelayer{front}
\pgfsetlayers{background,main,foreground,front}

\usepackage{multicol}
\usepackage{subcaption}
\newcommand{\pedge}[9]{
	
	\ifx\relax#6\relax
	\def\qoffs{0pt}
	\else
	\def\qoffs{#6}
	\fi
	
	\def\phedge{
		($#1+#5!\qoffs!-90:#2-#5$) -- 
		($#2+#1!\qoffs!-90:#3-#1$) -- 
		($#3+#2!\qoffs!-90:#4-#2$) -- 
		($#4+#3!\qoffs!-90:#5-#3$) -- 
		($#5+#4!\qoffs!-90:#1-#4$) -- cycle}

	\coordinate (12) at ($#1!\qoffs!90:#2$);
	\coordinate (15) at ($#1!\qoffs!-90:#5$);
	\coordinate (23) at ($#2!\qoffs!90:#3$);
	\coordinate (21) at ($#2!\qoffs!-90:#1$);
	\coordinate (34) at ($#3!\qoffs!90:#4$);
	\coordinate (32) at ($#3!\qoffs!-90:#2$);
	\coordinate (45) at ($#4!\qoffs!90:#5$);
	\coordinate (43) at ($#4!\qoffs!-90:#3$);
	\coordinate (51) at ($#5!\qoffs!90:#1$);
	\coordinate (54) at ($#5!\qoffs!-90:#4$);

	\def\nphedge{
		(15) let \p1=($(15)-#1$), \p2=($(12)-#1$) in 
		arc[start angle={atan2(\y1,\x1)}, delta angle={atan2(\y2,\x2)-atan2(\y1,\x1)-360*(atan2(\y2,\x2)-atan2(\y1,\x1)>0)}, x radius=\qoffs, y radius=\qoffs] --
		(21) let \p1=($(21)-#2$), \p2=($(23)-#2$) in 
		arc[start angle={atan2(\y1,\x1)}, delta angle={atan2(\y2,\x2)-atan2(\y1,\x1)-360*(atan2(\y2,\x2)-atan2(\y1,\x1)>0)}, x radius=\qoffs, y radius=\qoffs] --
		(32) let \p1=($(32)-#3$), \p2=($(34)-#3$) in 
		arc[start  angle={atan2(\y1,\x1)}, delta angle={atan2(\y2,\x2)-atan2(\y1,\x1)-360*(atan2(\y2,\x2)-atan2(\y1,\x1)>0)}, x radius=\qoffs, y radius=\qoffs] --
		(43) let \p1=($(43)-#4$), \p2=($(45)-#4$) in 
		arc[start angle={atan2(\y1,\x1)}, delta angle={atan2(\y2,\x2)-atan2(\y1,\x1)-360*(atan2(\y2,\x2)-atan2(\y1,\x1)>0)}, x radius=\qoffs, y radius=\qoffs] --
		(54) let \p1=($(54)-#5$), \p2=($(51)-#5$) in 
		arc[start angle={atan2(\y1,\x1)}, delta angle={atan2(\y2,\x2)-atan2(\y1,\x1)-360*(atan2(\y2,\x2)-atan2(\y1,\x1)>0)}, x radius=\qoffs, y radius=\qoffs] --
		cycle}

	\ifx\relax#7\relax
	\def\plwidth{1pt}
	\else
	\def\plwidth{#7}
	\fi
	
	\ifx\relax#9\relax
	\fill \nphedge;
	\else
	\fill[#9]\nphedge;
	\fi
	
	\ifx\relax#8\relax
	\draw[line width=\plwidth,rounded corners=\qoffs]\nphedge;
	\else
	\draw[line width=\plwidth,#8]\nphedge;
	\fi
}

\newcommand{\qedge}[7]{
	
	\ifx\relax#4\relax
	\def\qoffs{0pt}
	\else
	\def\qoffs{#4}
	\fi
	
	\def\qhedge{
		($#1+#3!\qoffs!-90:#2-#3$) --
		($#2+#1!\qoffs!-90:#3-#1$) --
		($#3+#2!\qoffs!-90:#1-#2$) -- cycle}

	\coordinate (12) at ($#1!\qoffs!90:#2$);
	\coordinate (13) at ($#1!\qoffs!-90:#3$);
	\coordinate (23) at ($#2!\qoffs!90:#3$);
	\coordinate (21) at ($#2!\qoffs!-90:#1$);
	\coordinate (31) at ($#3!\qoffs!90:#1$);
	\coordinate (32) at ($#3!\qoffs!-90:#2$);
	
	\def\nqhedge{
		(13) let \p1=($(13)-#1$), \p2=($(12)-#1$) in
		arc[start angle={atan2(\y1,\x1)}, delta angle={atan2(\y2,\x2)-atan2(\y1,\x1)-360*(atan2(\y2,\x2)-atan2(\y1,\x1)>0)}, x radius=\qoffs, y radius=\qoffs] --
		(21) let \p1=($(21)-#2$), \p2=($(23)-#2$) in
		arc[start angle={atan2(\y1,\x1)}, delta angle={atan2(\y2,\x2)-atan2(\y1,\x1)-360*(atan2(\y2,\x2)-atan2(\y1,\x1)>0)}, x radius=\qoffs, y radius=\qoffs] --
		(32) let \p1=($(32)-#3$), \p2=($(31)-#3$) in
		arc[start angle={atan2(\y1,\x1)}, delta angle={atan2(\y2,\x2)-atan2(\y1,\x1)-360*(atan2(\y2,\x2)-atan2(\y1,\x1)>0)}, x radius=\qoffs, y radius=\qoffs] --
		cycle}
	
	\ifx\relax#5\relax
	\def\qlwidth{1pt}
	\else
	\def\qlwidth{#5}
	\fi
	
	\ifx\relax#7\relax
	\fill \nqhedge;
	\else
	\fill[#7]\nqhedge;
	\fi
	
	\ifx\relax#6\relax
	\draw[line width=\qlwidth,rounded corners=\qoffs]\nqhedge;
	\else
	\draw[line width=\qlwidth,#6]\nqhedge;
	\fi
}

\newcommand{\redge}[8]{
	
	\ifx\relax#5\relax
	\def\qoffs{0pt}
	\else
	\def\qoffs{#5}
	\fi
	
	\def\rhedge{
		($#1+#4!\qoffs!-90:#2-#4$) -- 
		($#2+#1!\qoffs!-90:#3-#1$) -- 
		($#3+#2!\qoffs!-90:#4-#2$) -- 
		($#4+#3!\qoffs!-90:#1-#3$) -- cycle}

	\coordinate (12) at ($#1!\qoffs!90:#2$);
	\coordinate (14) at ($#1!\qoffs!-90:#4$);
	\coordinate (23) at ($#2!\qoffs!90:#3$);
	\coordinate (21) at ($#2!\qoffs!-90:#1$);
	\coordinate (34) at ($#3!\qoffs!90:#4$);
	\coordinate (32) at ($#3!\qoffs!-90:#2$);
	\coordinate (41) at ($#4!\qoffs!90:#1$);
	\coordinate (43) at ($#4!\qoffs!-90:#3$);
	
	\def\nrhedge{
		(14) let \p1=($(14)-#1$), \p2=($(12)-#1$) in 
		arc[start angle={atan2(\y1,\x1)}, delta angle={atan2(\y2,\x2)-atan2(\y1,\x1)-360*(atan2(\y2,\x2)-atan2(\y1,\x1)>0)}, x radius=\qoffs, y radius=\qoffs] --
		(21) let \p1=($(21)-#2$), \p2=($(23)-#2$) in 
		arc[start angle={atan2(\y1,\x1)}, delta angle={atan2(\y2,\x2)-atan2(\y1,\x1)-360*(atan2(\y2,\x2)-atan2(\y1,\x1)>0)}, x radius=\qoffs, y radius=\qoffs] --
		(32) let \p1=($(32)-#3$), \p2=($(34)-#3$) in 
		arc[start angle={atan2(\y1,\x1)}, delta angle={atan2(\y2,\x2)-atan2(\y1,\x1)-360*(atan2(\y2,\x2)-atan2(\y1,\x1)>0)}, x radius=\qoffs, y radius=\qoffs] --
		(43) let \p1=($(43)-#4$), \p2=($(41)-#4$) in 
		arc[start angle={atan2(\y1,\x1)}, delta angle={atan2(\y2,\x2)-atan2(\y1,\x1)-360*(atan2(\y2,\x2)-atan2(\y1,\x1)>0)}, x radius=\qoffs, y radius=\qoffs] --
		cycle}
	
	\ifx\relax#6\relax
	\def\rlwidth{1pt}
	\else
	\def\rlwidth{#6}
	\fi
	
	\ifx\relax#8\relax
	\fill \nrhedge;
	\else
	\fill[#8]\nrhedge;
	\fi
	
	\ifx\relax#7\relax
	\draw[line width=\rlwidth,rounded corners=\qoffs]\nrhedge;
	\else
	\draw[line width=\rlwidth,#7]\nrhedge;
	\fi
}

\let\epsilon=\varepsilon
\let\eps=\epsilon
\let\rho=\varrho
\let\theta=\vartheta

\newtheoremstyle{note}  {4pt}  {4pt}  {\sl}  {}  {\bfseries}  {.}  {.5em}          {}
\newtheoremstyle{introthms}  {3pt}  {3pt}  {\itshape}  {}  {\bfseries}  {.}  {.5em}          {\thmnote{#3}}
\newtheoremstyle{remark}  {2pt}  {2pt}  {\rm}  {}  {\bfseries}  {.}  {.3em}          {}

\theoremstyle{plain}
\newtheorem{theorem}{Theorem}[section]
\newtheorem{lemma}[theorem]{Lemma}
\newtheorem{prop}[theorem]{Proposition}

\theoremstyle{note}
\newtheorem{dfn}[theorem]{Definition}

\theoremstyle{remark}

\usepackage{lineno}
\newcommand*\patchAmsMathEnvironmentForLineno[1]{
	\expandafter\let\csname old#1\expandafter\endcsname\csname #1\endcsname
	\expandafter\let\csname oldend#1\expandafter\endcsname\csname end#1\endcsname
	\renewenvironment{#1}
	{\linenomath\csname old#1\endcsname}
	{\csname oldend#1\endcsname\endlinenomath}}
\newcommand*\patchBothAmsMathEnvironmentsForLineno[1]{
	\patchAmsMathEnvironmentForLineno{#1}
	\patchAmsMathEnvironmentForLineno{#1*}}
\AtBeginDocument{
	\patchBothAmsMathEnvironmentsForLineno{equation}
	\patchBothAmsMathEnvironmentsForLineno{align}
	\patchBothAmsMathEnvironmentsForLineno{flalign}
	\patchBothAmsMathEnvironmentsForLineno{alignat}
	\patchBothAmsMathEnvironmentsForLineno{gather}
	\patchBothAmsMathEnvironmentsForLineno{multline}
}

\def\ex{\text{\rm ex}}

\usepackage{scalerel}

\newsavebox\vdegbox
\savebox\vdegbox{\tikz{
		\draw[black,fill=black] (90:1) circle (.35);
		\draw[black,line width=0.10cm] (210:1) circle (.30);
		\draw[black,line width=0.10cm] (330:1) circle (.30);
		\draw[opacity=0] (0:1.2) circle (0.1);
	}}

\newsavebox\vvbox
\savebox\vvbox{\tikz{
		\draw[black,line width=0.10cm] (90:1) circle (.30);
		\draw[black,fill=black] (210:1) circle (.35);
		\draw[black,fill=black] (330:1) circle (.35);
		\draw[opacity=0] (0:1.2) circle (0.1);
	}}

\newsavebox\pdegbox
\savebox\pdegbox{\tikz{
		\draw[black,line width=0.10cm] (90:1) circle (.30);
		\draw[black,fill=black] (210:1) circle (.35);
		\draw[black,fill=black] (330:1) circle (.35);
		\draw[black,line width=0.28cm ] (210:1) -- (330:1);
		\draw[opacity=0] (0:1.2) circle (0.1);
	}}

\newsavebox\vvvbox
\savebox\vvvbox{\tikz{
		\draw[black,fill=black] (90:1) circle (.35);
		\draw[black,fill=black] (210:1) circle (.35);
		\draw[black,fill=black] (330:1) circle (.35);
		\draw[opacity=0] (0:1.2) circle (0.1);
	}}
\newcommand{\vvv}{\mathord{\scaleobj{1.2}{\scalerel*{\usebox{\vvvbox}}{x}}}}

\newsavebox\evbox
\savebox\evbox{\tikz{
		\draw[black,fill=black] (90:1) circle (.35);
		\draw[black,fill=black] (210:1) circle (.35);
		\draw[black,fill=black] (330:1) circle (.35);
		\draw[black,line width=0.28cm ] (210:1) -- (330:1);
		\draw[opacity=0] (0:1.2) circle (0.1);
	}}

\newsavebox\eebox
\savebox\eebox{\tikz{
		\draw[black,fill=black] (90:1) circle (.35);
		\draw[black,fill=black] (210:1) circle (.35);
		\draw[black,fill=black] (330:1) circle (.35);
		\draw[black,line width=0.28cm ] (90:1) -- (330:1);
		\draw[black,line width=0.28cm ] (90:1) -- (210:1);
		\draw[opacity=0] (0:1.2) circle (0.1);
	}}

\newsavebox\eeebox
\savebox\eeebox{\tikz{
		\draw[black,fill=black] (90:1) circle (.35);
		\draw[black,fill=black] (210:1) circle (.35);
		\draw[black,fill=black] (330:1) circle (.35);
		\draw[black,line width=0.28cm ] (90:1) -- (330:1);
		\draw[black,line width=0.28cm ] (90:1) -- (210:1);
		\draw[black,line width=0.28cm ] (210:1) -- (330:1);
		\draw[opacity=0] (0:1.2) circle (0.1);
	}}

\makeatletter
\newcommand{\overrighharpoonup}[1]{\ThisStyle{%
		\vbox {\m@th\ialign{##\crcr
				\rightharpoonupfill \crcr
				\noalign{\kern-\p@\nointerlineskip}
				$\hfil\SavedStyle#1\hfil$\crcr}}}}

\def\rightharpoonupfill{%
	$\SavedStyle\m@th\mkern+0.8mu\cleaders\hbox{$\shortbar\mkern-4mu$}\hfill\rightharpoonuptip\mkern+0.8mu$}

\def\rightharpoonuptip{%
	\raisebox{\z@}[2pt][1pt]{\scalebox{0.55}{$\SavedStyle\rightharpoonup$}}}

\def\shortbar{%
	\smash{\scalebox{0.55}{$\SavedStyle\relbar$}}}
\makeatother

\makeatletter
\newcommand{\overlefharpoonup}[1]{\ThisStyle{%
		\vbox {\m@th\ialign{##\crcr
				\leftharpoonupfill \crcr
				\noalign{\kern-\p@\nointerlineskip}
				$\hfil\SavedStyle#1\hfil$\crcr}}}}

\def\leftharpoonupfill{%
	$\SavedStyle\m@th\mkern+0.8mu\cleaders\hbox{$\shortbar\mkern-4mu$}\hfill\leftharpoonuptip\mkern+0.8mu$}

\def\leftharpoonuptip{%
	\raisebox{\z@}[2pt][1pt]{\scalebox{0.55}{$\SavedStyle\leftharpoonup$}}}

\makeatletter
\newsavebox\myboxA
\newsavebox\myboxB
\newlength\mylenA

\newcommand*\xoverline[2][0.75]{%
	\sbox{\myboxA}{$\m@th#2$}%
	\setbox\myboxB\null
	\ht\myboxB=\ht\myboxA%
	\dp\myboxB=\dp\myboxA%
	\wd\myboxB=#1\wd\myboxA
	\sbox\myboxB{$\m@th\overline{\copy\myboxB}$}
	\setlength\mylenA{\the\wd\myboxA}
	\addtolength\mylenA{-\the\wd\myboxB}%
	\ifdim\wd\myboxB<\wd\myboxA%
	\rlap{\hskip 0.5\mylenA\usebox\myboxB}{\usebox\myboxA}%
	\else
	\hskip -0.5\mylenA\rlap{\usebox\myboxA}{\hskip 0.5\mylenA\usebox\myboxB}%
	\fi}
\makeatother

\begin{document}
	
	\title[Infinitely many accumulation points of codegree Tur\'an densities]
	{Infinitely many accumulation points of codegree Tur\'an densities}

\author[Heng Li]{Heng Li}
    \address{School of Mathematics, Shandong University, Jinan, China}
    \email{heng.li@sdu.edu.cn}

\author[Weichan Liu]{Weichan Liu}
	\address{School of Mathematics, Shandong University, Jinan, China}
    \email{wcliu@sdu.edu.cn}

\author[Bjarne Sch\"{u}lke]{Bjarne Sch\"{u}lke}
	\address{Extremal Combinatorics and Probability Group, Institute for Basic Science, Daejeon, South Korea}
	\email{schuelke@ibs.re.kr}
 
 \author[Wanting Sun]{Wanting Sun}
	\address{Data Science Institute, Shandong University, Jinan, China}
	\email{wtsun@sdu.edu.cn}


\begin{abstract}
The codegree Tur\'an density~$\gamma(F)$ of a~$k$-graph~$F$ is the smallest~$\gamma\in[0,1)$ such that every~$k$-graph~$H$ with~$\delta_{k-1}(H)\geq(\gamma+o(1))\vert V(H)\vert$ contains a copy of~$F$.
We prove that for all~$k,r\in\mathds{N}$ with~$k\geq3$,~$\frac{r-1}{r}$ is an accumulation point of~$\Gamma^{(k)}=\left\{\gamma(F):F\text{ is a }k\text{-graph}\right\}$.
This makes progress on a problem posed by Mubayi and Zhao.
\end{abstract}

\maketitle
	
\section{Introduction}\label{SEC:Introduction}

A $k$-uniform hypergraph (or~$k$-graph) $H$ consists of a vertex set~$V(H)$ together with a set of edges~$E(H)\subseteq V(H)^{(k)}=\left\{S\subseteq V(H):\vert S\vert =k\right\}$.
	Given a~$k$-graph~$F$ and~$n\in\mathds{N}$, the Tur\'an number of~$n$ and~$F$, denoted by~$\ex(n,F)$, is the maximum number of edges an~$n$-vertex~$k$-graph can have without containing a copy of~$F$. 
    Since the main interest lies in the asymptotics, the \emph{Tur\'an density}~$\pi(F)$ of a~$k$-graph~$F$ is defined as
    $$\pi(F) = \displaystyle\lim_{n \to \infty} \frac{\ex(n,F)}{\binom{n}{k}}\,.$$
	Determining the value of~$\pi(F)$ for~$k$-graphs (with~$k\geq 3$) is one of the central open problems in combinatorics. 
	In particular, the problem of determining the Tur\'an density of the complete~$3$-graph on four vertices, i.e.,~$\pi(K_4^{(3)})$, was asked by Tur\'an in 1941~\cite{T:41} and Erd\H{o}s~\cite{E:77} offered {\$500} for determining any~$\pi(K_{\ell}^{(k)})$ with~$\ell > k\geq 3$ and~{\$1000} for determining all~$\pi(K_{\ell}^{(k)})$ with~$\ell >  k\geq 3$.
	Despite receiving a lot of attention (see, for instance, the surveys on the topic~\cites{F:91,K:11,S:95}), this problem, and even the seemingly simpler problem of determining~$\pi(K_4^{(3)-})$, where~$K_4^{(3)-}$ is the~$K_4^{(3)}$ minus one edge, remain open.
	Several variations of this type of problem have been considered, the most prominent ones being the uniform Tur\'an density and the codegree Tur\'an density, see for instance,~\cites{BCL:21,ES:82,MZ:07,R:20}.
	
    Here we are concerned with the codegree Tur\'an density, which asks how large the \emph{minimum codegree} of an~$F$-free~$k$-graph can be.
	Given a~$k$-graph~$H=(V,E)$ and~$S\subseteq V$, the degree~$d(S)$ of~$S$ is the number of edges containing~$S$, i.e., $d(S)=\vert\{e\in E:S\subseteq e\}\vert$.
	The \emph{minimum codegree} of~$H$ is defined as~$\delta(H)=\min_{x\in V^{(k-1)}}d(x)$\footnote{When considering different types of degree for~$k$-graphs, this is usually denoted as~$\delta_{k-1}(H)$. Here we suppress the index to simplify the notation.}.
	Given a~$k$-graph~$F$ and~$n\in\mathds{N}$, Mubayi and Zhao \cite{MZ:07} introduced the \emph{codegree Tur\'an number}~$\ex_{\text{co}}(n,F)$ of~$n$ and~$F$ as the maximum~$d$ such that there is an~$F$-free $k$-graph~$H$ on~$n$ vertices with~$\delta(H)\geq d$. 
	Moreover, they defined the \emph{codegree Tur\'an density}~$\gamma(F)$ of $F$ as
	$$\gamma(F) = \lim_{n\to\infty} \frac{\ex_{\text{co}}(n,F)}{n}\,$$
    and proved that this limit always exists. 
    It is not hard to see that~$\gamma(F) \leq \pi(F)\,.$
    The codegree Tur\'an density of a family~$\mathcal{F}$ of~$k$-graphs is defined analogously.
    Similarly as for the Tur\'an density, determining the exact codegree Tur\'an density of a given hypergraph is in general very difficult and so it is only known for very few hypergraphs (see the table in~\cite{BCL:21}).
    In order to study the behaviour of possible codegree Tur\'an densities, Mubayi and Zhao~\cite{MZ:07} considered the sets
    \begin{align*}
        \Gamma^{(k)}&:=\left\{\gamma(F) \colon F \text{ is a~$k$-graph}\right\}\,,\\
        \Gamma^{(k)}_{\text{fin}}&:=\{\gamma(\mathcal{F}) \colon \mathcal{F} \text{ is a finite family of~$k$-graphs}\}\,,\text{ and }\\
        \Gamma^{(k)}_{\infty}&:=\{\gamma(\mathcal{F}) \colon \mathcal{F} \text{ is a family of~$k$-graphs}\}\,.
    \end{align*}
    Note that~$\Gamma^{(k)}\subseteq\Gamma^{(k)}_{\text{fin}}\subseteq\Gamma^{(k)}_{\infty}\subseteq[0,1)$.
    The investigation of the analogous sets~$\Pi^{(k)}$,~$\Pi^{(k)}_{\text{fin}}$, and~$\Pi^{(k)}_{\infty}$ for the classical Tur\'an density has a rich history that is intimately tied to the famous Erd\H{o}s jumping conjecture (see, for instance,~\cites{BT:11,CS:25,FR:84,P:14,S:23}).
    More recently, there have also been some results on the respective sets~$\Pi_{\vvv}$,~$\Pi_{\vvv,\text{fin}}$, and~$\Pi_{\vvv,\infty}$ for the uniform Tur\'an density of~$3$-graphs~\cites{KSS:24,L:24,RRS:18}.
    As Pikhurko emphasises in~\cite{P:14}, the case in which one single~$k$-graph is forbidden is arguably the most interesting one, meaning the study of~$\Pi^{(k)}$ and~$\Gamma^{(k)}$. 
    However, as this case is in general notoriously difficult, many researchers have so far directed their attention to the set of possible Tur\'an densities of families of hypergraphs.
    For each of~$\Gamma^{(k)}$,~$\Pi^{(k)}$, and~$\Pi_{\vvv}$, it has only been shown very recently that this set has an accumulation point (in~$[0,1)$) in~\cite{PS:23},~\cite{CS:25}, and~\cite{L:24}, respectively.
    In stark contrast, a classical theorem by Erd\H{o}s, Stone, and Simonovits~\cite{ES:46,ES:66} states that~$\Pi^{(2)}=\{\frac{r-1}{r}:\:r\in\mathds{N}\}$.
    This indicates the difficulty of Tur\'an-type problems for hypergraphs.

    Mubayi and Zhao~\cite{MZ:07} proved that~$\Gamma^{(k)}_{\text{fin}}$ is dense in~$[0,1]$ and they propose the problem of describing~$\Gamma^{(k)}$. 
    In particular, they ask whether~$\Gamma^{(k)}$ is dense in~$[0,1]$.
    In~\cite{PS:23} a first step towards this was taken by proving the existence of an accumulation point of~$\Gamma^{(k)}\subseteq[0,1)$.
    \begin{theorem}[\cite{PS:23}]\label{thm:zero}
        Zero is an accumulation point of~$\Gamma^{(k)}$ for every integer~$k\geq3$.
    \end{theorem}

    Here we make further progress on the aforementioned problem by showing that~$\Gamma^{(k)}$ has infinitely many accumulation points.
    \begin{theorem}\label{thm:main}
    For all~$k, r\in\mathds{N}$ with $k\geq 3$,
    ~$\frac{r-1}{r}$ is an accumulation point of~$\Gamma^{(k)}$.
    \end{theorem}

    We remark that Keevash and Zhao~\cite{KZ:07} showed that~$\frac{r-1}{r}\in\Gamma^{(k)}$ for every~$r\in\mathds{N}$ and~$k\in\mathds{N}$.
    Since this paper first appeared online, Conlon and Sch\"ulke~\cite{CS:25} showed that~$\Pi^{(k)}$ has infinitely many accumulation points for~$k\geq3$.
    However, the proof ideas are very different and the accumulation points there are accumulation points ``from below'' whereas the accumulation points in Theorem~\ref{thm:main} are ``from above''.

\section{Preliminaries}
Given a~$k$-graph~$H=(V,E)$ and~$S\in V^{(k-1)}$, the (co)neighbourhood of~$S$ is~$N(S)=\{v\in V:\{v\}\cup S\in E\}$.
If~$S=\{v_1,\dots,v_{k-1}\}$ for some vertices~$v_1,\dots,v_{k-1}\in V$, then we simply write~$N(v_1\cdots v_{k-1})$ instead of~$N(\{v_1,\dots ,v_{k-1}\})$.

The following proposition states that the supersaturation phenomenon occurs also for the codegree Tur\'an density.
It was effectively proved by Mubayi and Zhao~\cite{MZ:07} (for this formulation and a proof see, for instance,~\cite{PS:23}).
Given~$t\in\mathds{N}$ and a~$k$-graph~$F$, let~$F(t)$ be the~$t$-blow-up of~$F$, i.e., the~$k$-graph
obtained from~$F$ by replacing every vertex by~$t$ copies of itself.
\begin{prop}\label{prop:supersaturation}
    Let~$t,k,c\in \mathds N$ with~$k\geq 2$ and let~$\mathcal F=\{F_1, \dots, F_t\}$ be a finite family of~$k$-graphs with~$\vert V(F_i)\vert=f_i$ for all $i\in [t]$.
    For every~$\eps>0$, there exists a~$\zeta>0$ such that for sufficiently large $n\in \mathds N$ the following holds.
    Every $n$-vertex~$k$-graph $H$ with~$\delta(H)\geq (\gamma(\mathcal F)+\eps)n$
    contains~$\zeta \binom{n}{f_i}$ copies of~$F_i$ for some~$i\in[t]$.
    Consequently, $H$ contains a copy of~$F_i(c)$.
\end{prop}

Given two~$k$-graphs~$F$ and~$G$, a \emph{homomorphism} from~$F$ to~$G$ is a map~$\varphi :\: V(F) \to V(G)$ such that~$\varphi(e) \in E(G)$ for all $e \in E(F)$. It easily follows from Proposition~\ref{prop:supersaturation} that if there is a homomorphism from~$F$ to~$G$, then~$\gamma(F)\leq\gamma(G)$.

In~\cite{PS:23}, Piga and Sch\"ulke used the following~$k$-graphs to prove that zero is an accumulation point for~$\Gamma^{(k)}$ for every integer~$k\geq3$.

\begin{dfn}\label{dfn:Zell}
    For integers~$\ell \geq k\geq 2$, the \emph{$k$-uniform zycle of length $\ell$} is the~$k$-graph~$Z_\ell^{(k)}$ given by 
    \begin{align*}
        V(Z_\ell^{(k)})=&\left\{v_i^j\colon i\in [\ell], j\in [k-1] \right\} \text{, and}\\
        E(Z_\ell^{(k)})=&\left\{v_i^1v_i^2\cdots v_i^{k-1} v_{i+1}^j \colon i\in [\ell], j\in[k-1]\right\}\,,
    \end{align*}    
    where the indices are considered modulo~$\ell$.
\end{dfn}

Let us remark that zycles have been considered independently before, see for instance~\cite{DJ:14}.
The following lemma provides an upper bound on the codegree Tur\'an density of some zycles.

\begin{lemma}[Lemma 2.1 in~\cite{PS:23}]\label{lem:zyclesupper}
    Let~$k\geq 3$.
    For every~$d\in (0,1]$, there is some~$\ell\in \mathds N$ such that
    $$\gamma \big(Z^{(k)}_\ell\big) \leq d\,.$$
\end{lemma}

For a lower bound on the codegree Tur\'an density of zycles, Piga and Sch\"ulke~\cite{PS:23} introduced the following~$k$-graphs.

\begin{dfn}\label{dfn:alggraph}
    Let~$n,p,k\in \mathds N$ be such that $p$ is a prime,~$k \geq 2$ and~$p\mid n$.
    We define the $n$-vertex~$k$-graph~$\mathds F_p^{(k)}(n)$ as follows.
    The vertex set consists of~$p$ disjoint sets of size~$\frac{n}{p}$ each, i.e., $V(\mathds F_p^{(k)}(n)) =  V_0 \dcup \dots \dcup V_{p-1}$ with~$\vert V_i\vert=\frac{n}{p}$ for all~$i\in [p]$.
    Given a vertex~$v\in V(\mathds F_p^{(k)}(n))$, we write~$\mathfrak f(v)=i$ if and only if~$v\in V_i$ for~$i\in \{0,1,\dots, p-1\}$.
    We define the edge set of~$\mathds F_p^{(k)}(n)$~by
    $$v_1\cdots v_k \in E(\mathds F_p^{(k)}(n)) \Leftrightarrow 
    \begin{cases}
        \mathfrak f(v_1)+\dots+ \mathfrak f(v_k) \equiv 0 \bmod p \text{ and } \mathfrak f(v_i) \neq 0 \text{ for some }i\in [k] \text{, or}\\
        \mathfrak f(v_{\sigma(1)})=\dots=\mathfrak f(v_{\sigma(k-1)})=0 \text{ and } \mathfrak f(v_{\sigma(k)})=1 \text{ for some }\sigma\in S_k\,.
    \end{cases}$$
\end{dfn}

They proved the following properties (see Lemma 2.4 and Display (2.4) in~\cite{PS:23}).

\begin{lemma}\label{lem:Fpprop}
    For every integer~$k\geq3$ and every~$\ell\in\mathds{N}$, there is some~$p_0$ such that for every prime~$p\geq p_0$ and every~$n\in\mathds{N}$ with~$p\mid n$, we have~$Z^{(k)}_{\ell}\not\subseteq \mathds{F}_p^{(k)}(n)$ and~$\delta\left(\mathds{F}_p^{(k)}(n)\right)=\frac{n}{p}$.
\end{lemma}

\section{Proof of Theorem~\ref{thm:main}}
\noindent{\bf Proof of Theorem~\ref{thm:main}.}
    Let~$k\geq3$ be an integer.
    For~$\eta>0$ and~$r,n\in\mathds{N}$ let~$H^r(\eta,n)$ be the~$k$-graph on~$n$ vertices with vertex set\footnote{For ease of notation, we suppress the dependencies of~$V$ and~$V_i$ on the parameters~$r$,~$\eta$, and~$n$.}~$V=V_1\dcup\dots\dcup V_r$, where~$\vert V_i\vert=\big\lfloor\frac{(1-\eta)n}{r}\big\rfloor$ for every~$i\in[r-1]$ (and~$\vert V_r\vert=n-\sum_{i\in[r-1]}\vert V_i\vert$), and edge set $$\left\{e\in V^{(k)}:\:e\not\subseteq V_i\text{ for all $i\in[r]$}\right\}\,.$$
    Further, if in addition~$p$ is a prime with~$p\mid \vert V_r\vert$, let~$H^r_{\mathds{F}_p}(\eta,n)$ be the~$k$-graph obtained from~$H^r(\eta,n)$ by placing a copy of~$\mathds{F}_p^{(k)}(\vert V_r\vert)$ inside~$V_r$.
    Note that~$H^1_{\mathds{F}_p}(\eta,n)=\mathds{F}_p^{(k)}(n)$.
    
    We will perform an induction on~$r$ to construct, for all integers~$\ell\geq2$,~$k$-graphs~$G_{\ell}^r$ with the following properties.
    \begin{enumerate}
        \item\label{it:upperbound} For every~$r\in\mathds{N}$ and every~$\varepsilon>0$, there is some~$\ell\in\mathds{N}$ such that~$\gamma(G_{m\cdot\ell}^r)\leq\frac{r-1}{r}+\varepsilon$ holds for every~$m\in\mathds{N}$.
        \item\label{it:lowerbound} For every~$r\in\mathds{N}$ and every integer~$\ell\geq 2$, there is some~$p_0$ such that~$G^r_{\ell}\not\subseteq H^r_{\mathds{F}_p}(\eta,n)$ for every prime~$p\geq p_0$,~$\eta>0$, and~$n\in\mathds{N}$ with~$p\mid n-(r-1)\cdot\big\lfloor\frac{(1-\eta)n}{r}\big\rfloor$.
    \end{enumerate}
    Before constructing the~$k$-graphs~$G_{\ell}^r$, let us observe that Theorem~\ref{thm:main} will follow from the existence of~$k$-graphs with these properties.
    Indeed, for this we only use property~\eqref{it:upperbound} with~$m=1$ and property~\eqref{it:lowerbound}; the stronger version of property~\eqref{it:upperbound} is used for the induction.
    
    Let~$\varepsilon>0$ and~$r\in\mathds{N}$ be given.
    We intend to show that there is some integer~$\ell$ such that~$\frac{r-1}{r}<\gamma(G_{\ell}^r)\leq\frac{r-1}{r}+\varepsilon$.
    By property~\eqref{it:upperbound}, there exists some~$\ell\in\mathds{N}$ such that~$\gamma(G_{\ell}^r)\leq\frac{r-1}{r}+\varepsilon$.
    By property~\eqref{it:lowerbound}, there is some~$p_0$ such that~$G^r_{\ell}\not\subseteq H^r_{\mathds{F}_p}(\eta,n)$ for every prime~$p\geq p_0$,~$\eta>0$, and~$n\in\mathds{N}$ with~$p\mid n-(r-1)\cdot\big\lfloor\frac{(1-\eta)n}{r}\big\rfloor$.
    Therefore, it remains to show that there are a prime~$p\geq p_0$,~$\eta>0$, and infinitely many~$n\in\mathds{N}$ with~$p\mid n-(r-1)\cdot\big\lfloor\frac{(1-\eta)n}{r}\big\rfloor$ such that~$\delta\big(H^r_{\mathds{F}_p}(\eta,n)\big)>\big(\frac{r-1}{r}+\frac{\eta}{r}\big)n$.
    For this, fix a prime~$p\geq p_0$ and choose~$\eta\in\mathds{Q}$ and any of the infinitely many choices for~$n\in\mathds{N}$ such that~$p\mid n-(r-1)\cdot\big\lfloor\frac{(1-\eta)n}{r}\big\rfloor$ and $$r,\ell,k, p\ll\frac{1}{\eta}\ll n\,.$$
    First, notice that for any~$S\in V^{(k-1)}$ with~$S\not\subseteq V_i$ for all~$i\in[r]$, we have~$d(S)=n-(k-1)\geq(1-\eta)n$.
    Second, for any~$S\in V^{(k-1)}$ with~$S\subseteq V_i$ for some~$i\in[r-1]$, we have~$d(S)\geq\frac{r-1}{r}n+\frac{\eta}{r}n$.
    Using Lemma~\ref{lem:Fpprop} one can see that for~$S\in V_r^{(k-1)}$, we have $$d(S)\geq(1-2\eta)\frac{r-1}{r}n+\frac{n}{rp}\geq\frac{r-1}{r}n+\frac{n}{2rp}\,.$$
    Thus,~$\delta\big(H^r_{\mathds{F}_p}(\eta,n)\big)\geq\left(\frac{r-1}{r}+\frac{\eta}{r}\right)n$.
    Hence, properties~\eqref{it:upperbound} and~\eqref{it:lowerbound} indeed imply the statement of the theorem.

    Next, we construct the~$k$-graphs~$G_{\ell}^r$ as claimed.
    For~$r=1$ and an integer~$\ell\geq2$, let~$G_{\ell}^1=Z_{\ell}^{(k)}$.
    As observed in~\cite{PS:23}, for every~$m\in\mathds{N}$ there is a sufficiently large blow-up of~$Z_{\ell}^{(k)}$ that contains a copy of~$Z_{m\cdot\ell}^{(k)}$.
    Hence, Proposition~\ref{prop:supersaturation} ensures that we have~$\gamma\big(Z_{m\cdot\ell}^{(k)}\big)\leq\gamma\big(Z_{\ell}^{(k)}\big)$ for every~$m\in\mathds{N}$.
    Then property~\eqref{it:upperbound} follows from Lemma~\ref{lem:zyclesupper}.
    Moreover, since~$H^1_{\mathds{F}_p}(\eta,n)$ is simply~$\mathds{F}^{(k)}_p(n)$, property~\eqref{it:lowerbound} follows from Lemma~\ref{lem:Fpprop}.

    Now assume that for some integer~$r>1$ and all integers~$\ell\geq2$, we have constructed~$k$-graphs~$G_{\ell}^{r-1}$ that satisfy the properties above.
    Consider~$G_{\ell}^{r-1}(k-1)$ and denote by~$x_v^1,x_v^2,\ldots,x_v^{k-1}$ the~$k-1$ copies of a vertex~$v\in V(G_{\ell}^{r-1})$.
    Roughly speaking, we obtain~$G_{\ell}^r$ from~$G_{\ell}^{r-1}(k-1)$ by adding, for each collection~$S$ of~$r-1$ $(k-1)$-sets of the form~$\{x_v^1, \dots, x_v^{k-1}\}$, a copy of~$Z_{\ell}^{(k)}$ on a new vertex set~$V_S$, as well as all edges from any~$(k-1)$-set in~$S$ to every vertex of~$V_S$.
    More formally, add to~$G_{\ell}^{r-1}(k-1)$ copies of~$Z_{\ell}^{(k)}$ on vertex sets~$V_S$ for every~$S\in V(G_{\ell}^{r-1})^{(r-1)}$ such that~$V_S\cap V_{S'}=\emptyset$ and~$V_S\cap V(G_{\ell}^{r-1}(k-1))=\emptyset$ for all distinct~$S,S'\in V(G_{\ell}^{r-1})^{(r-1)}$.
    Further, for every~$S\in V(G_{\ell}^{r-1})^{(r-1)}$ and~$v\in S$, add all edges of the form~$x_v^1x_v^2\cdots x_v^{k-1}y$ with~$y\in V_S$.
    Call the resulting~$k$-graph~$G_{\ell}^r$.
    See Figure \ref{FIGHRE:G_L^R} for an illustration of this construction.

\begin{figure}
    \centering

\tikzset{every picture/.style={line width=0.75pt}} 

\tikzset{every picture/.style={line width=0.75pt}} 

\tikzset{every picture/.style={line width=0.75pt}} 

\tikzset{every picture/.style={line width=0.75pt}} 

\tikzset{every picture/.style={line width=0.75pt}} 

\tikzset{every picture/.style={line width=0.75pt}} 

\begin{tikzpicture}[x=0.75pt,y=0.75pt,yscale=-0.95,xscale=0.95]

\draw  [draw opacity=0][fill={rgb, 255:red, 155; green, 155; blue, 155 }  ,fill opacity=0.42 ][dash pattern={on 0.84pt off 2.51pt}][line width=0.75]  (351.39,406) .. controls (354.05,406) and (356.2,410.05) .. (356.2,415.05) .. controls (356.2,420.05) and (354.05,424.1) .. (351.39,424.1) .. controls (348.73,424.1) and (346.57,420.05) .. (346.57,415.05) .. controls (346.57,410.05) and (348.73,406) .. (351.39,406) -- cycle ;
\draw [draw opacity=0][fill={rgb, 255:red, 208; green, 2; blue, 27 }  ,fill opacity=0.5 ]   (357.3,415.66) .. controls (372.81,434.28) and (393.87,448.73) .. (389.43,447.17) .. controls (384.98,445.61) and (361.12,440.15) .. (357.3,445.66) .. controls (353.47,451.17) and (366.14,434.06) .. (357.3,415.66) -- cycle ;
\draw [draw opacity=0][fill={rgb, 255:red, 208; green, 2; blue, 27 }  ,fill opacity=0.5 ]   (77.43,420.77) .. controls (97.96,419.79) and (122.61,403.1) .. (120.15,406.29) .. controls (117.68,409.48) and (89.61,444.67) .. (91.1,445.21) .. controls (92.59,445.74) and (95.08,434.16) .. (77.43,420.77) -- cycle ;
\draw [draw opacity=0][fill={rgb, 255:red, 208; green, 2; blue, 27 }  ,fill opacity=0.5 ]   (77.43,420.77) .. controls (99.66,430) and (122.11,434.88) .. (122.45,433.9) .. controls (122.79,432.91) and (91.27,444.95) .. (91.1,445.21) .. controls (90.94,445.46) and (95.08,434.16) .. (77.43,420.77) -- cycle ;
\draw  [fill={rgb, 255:red, 0; green, 0; blue, 0 }  ,fill opacity=1 ] (76.45,420.97) .. controls (76.18,420.48) and (76.37,419.87) .. (76.85,419.61) .. controls (77.34,419.35) and (77.95,419.53) .. (78.21,420.02) .. controls (78.47,420.51) and (78.29,421.11) .. (77.8,421.37) .. controls (77.31,421.64) and (76.71,421.45) .. (76.45,420.97) -- cycle ;
\draw  [fill={rgb, 255:red, 0; green, 0; blue, 0 }  ,fill opacity=1 ] (90.12,445.4) .. controls (89.86,444.92) and (90.04,444.31) .. (90.53,444.05) .. controls (91.02,443.79) and (91.62,443.97) .. (91.88,444.46) .. controls (92.15,444.94) and (91.96,445.55) .. (91.48,445.81) .. controls (90.99,446.07) and (90.38,445.89) .. (90.12,445.4) -- cycle ;
\draw [draw opacity=0][fill={rgb, 255:red, 208; green, 2; blue, 27 }  ,fill opacity=0.5 ]   (162.05,411.67) .. controls (176.36,426.44) and (203.37,429.43) .. (199.35,429.69) .. controls (195.32,429.95) and (152.16,436.62) .. (152.74,438.08) .. controls (153.32,439.55) and (163.67,433.77) .. (162.05,411.67) -- cycle ;
\draw [draw opacity=0][fill={rgb, 255:red, 208; green, 2; blue, 27 }  ,fill opacity=0.5 ]   (162.05,411.67) .. controls (169.81,434.46) and (180.98,454.54) .. (181.94,454.14) .. controls (182.9,453.75) and (153.04,438.04) .. (152.74,438.08) .. controls (152.44,438.13) and (163.67,433.77) .. (162.05,411.67) -- cycle ;
\draw [draw opacity=0][fill={rgb, 255:red, 208; green, 2; blue, 27 }  ,fill opacity=0.5 ]   (120.15,406.29) .. controls (138.71,415.13) and (165.06,408.5) .. (161.38,410.15) .. controls (157.71,411.8) and (119.59,433.12) .. (120.65,434.29) .. controls (121.71,435.46) and (129.38,426.44) .. (120.15,406.29) -- cycle ;
\draw [draw opacity=0][fill={rgb, 255:red, 208; green, 2; blue, 27 }  ,fill opacity=0.5 ]   (120.15,406.29) .. controls (135.38,424.94) and (152.85,439.85) .. (153.62,439.14) .. controls (154.38,438.44) and (120.91,434.15) .. (120.65,434.29) .. controls (120.38,434.44) and (129.38,426.44) .. (120.15,406.29) -- cycle ;
\draw  [fill={rgb, 255:red, 0; green, 0; blue, 0 }  ,fill opacity=1 ] (197.93,429.77) .. controls (198.3,429.36) and (198.94,429.32) .. (199.35,429.69) .. controls (199.76,430.06) and (199.79,430.69) .. (199.42,431.1) .. controls (199.06,431.51) and (198.42,431.55) .. (198.01,431.18) .. controls (197.6,430.81) and (197.57,430.18) .. (197.93,429.77) -- cycle ;
\draw  [fill={rgb, 255:red, 0; green, 0; blue, 0 }  ,fill opacity=1 ] (180.45,452.81) .. controls (180.82,452.4) and (181.45,452.36) .. (181.86,452.73) .. controls (182.28,453.1) and (182.31,453.73) .. (181.94,454.14) .. controls (181.57,454.56) and (180.94,454.59) .. (180.53,454.22) .. controls (180.12,453.85) and (180.08,453.22) .. (180.45,452.81) -- cycle ;
\draw  [fill={rgb, 255:red, 0; green, 0; blue, 0 }  ,fill opacity=1 ] (119.19,406) .. controls (119.19,405.45) and (119.64,405) .. (120.19,405) .. controls (120.74,405) and (121.19,405.45) .. (121.19,406) .. controls (121.19,406.55) and (120.74,407) .. (120.19,407) .. controls (119.64,407) and (119.19,406.55) .. (119.19,406) -- cycle ;
\draw  [fill={rgb, 255:red, 0; green, 0; blue, 0 }  ,fill opacity=1 ] (119.69,434) .. controls (119.69,433.45) and (120.14,433) .. (120.69,433) .. controls (121.24,433) and (121.69,433.45) .. (121.69,434) .. controls (121.69,434.55) and (121.24,435) .. (120.69,435) .. controls (120.14,435) and (119.69,434.55) .. (119.69,434) -- cycle ;
\draw  [fill={rgb, 255:red, 0; green, 0; blue, 0 }  ,fill opacity=1 ] (160.43,409.86) .. controls (160.59,409.33) and (161.15,409.03) .. (161.68,409.19) .. controls (162.2,409.35) and (162.5,409.91) .. (162.34,410.44) .. controls (162.18,410.97) and (161.62,411.27) .. (161.09,411.11) .. controls (160.57,410.94) and (160.27,410.39) .. (160.43,409.86) -- cycle ;
\draw  [fill={rgb, 255:red, 0; green, 0; blue, 0 }  ,fill opacity=1 ] (151.7,438.56) .. controls (151.86,438.03) and (152.42,437.73) .. (152.95,437.89) .. controls (153.48,438.06) and (153.78,438.61) .. (153.62,439.14) .. controls (153.46,439.67) and (152.9,439.97) .. (152.37,439.81) .. controls (151.84,439.65) and (151.54,439.09) .. (151.7,438.56) -- cycle ;
\draw  [fill={rgb, 255:red, 0; green, 0; blue, 0 }  ,fill opacity=1 ] (350.37,410.36) .. controls (350.16,409.84) and (350.42,409.26) .. (350.93,409.06) .. controls (351.45,408.86) and (352.03,409.11) .. (352.23,409.63) .. controls (352.43,410.14) and (352.18,410.72) .. (351.66,410.92) .. controls (351.15,411.13) and (350.57,410.87) .. (350.37,410.36) -- cycle ;
\draw  [fill={rgb, 255:red, 0; green, 0; blue, 0 }  ,fill opacity=1 ] (350.37,420.36) .. controls (350.16,419.84) and (350.42,419.26) .. (350.93,419.06) .. controls (351.45,418.86) and (352.03,419.11) .. (352.23,419.63) .. controls (352.43,420.14) and (352.18,420.72) .. (351.66,420.92) .. controls (351.15,421.13) and (350.57,420.87) .. (350.37,420.36) -- cycle ;
\draw  [fill={rgb, 255:red, 0; green, 0; blue, 0 }  ,fill opacity=1 ] (350.37,440.36) .. controls (350.16,439.84) and (350.42,439.26) .. (350.93,439.06) .. controls (351.45,438.86) and (352.03,439.11) .. (352.23,439.63) .. controls (352.43,440.14) and (352.18,440.72) .. (351.66,440.92) .. controls (351.15,441.13) and (350.57,440.87) .. (350.37,440.36) -- cycle ;
\draw  [fill={rgb, 255:red, 0; green, 0; blue, 0 }  ,fill opacity=1 ] (350.37,450.36) .. controls (350.16,449.84) and (350.42,449.26) .. (350.93,449.06) .. controls (351.45,448.86) and (352.03,449.11) .. (352.23,449.63) .. controls (352.43,450.14) and (352.18,450.72) .. (351.66,450.92) .. controls (351.15,451.13) and (350.57,450.87) .. (350.37,450.36) -- cycle ;
\draw  [dash pattern={on 4.5pt off 4.5pt}]  (74.04,448.13) .. controls (62.04,548.13) and (201.37,559.46) .. (198.37,452.8) ;
\draw [draw opacity=0][fill={rgb, 255:red, 208; green, 2; blue, 27 }  ,fill opacity=0.5 ]   (357.3,415.66) .. controls (375.2,423.39) and (404.54,423.17) .. (399.43,423.17) .. controls (394.32,423.17) and (361.12,440.15) .. (357.3,445.66) .. controls (353.47,451.17) and (366.14,434.06) .. (357.3,415.66) -- cycle ;
\draw  [draw opacity=0][fill={rgb, 255:red, 155; green, 155; blue, 155 }  ,fill opacity=0.42 ][dash pattern={on 0.84pt off 2.51pt}][line width=0.75]  (351.39,435.63) .. controls (354.05,435.63) and (356.2,439.68) .. (356.2,444.68) .. controls (356.2,449.68) and (354.05,453.73) .. (351.39,453.73) .. controls (348.73,453.73) and (346.57,449.68) .. (346.57,444.68) .. controls (346.57,439.68) and (348.73,435.63) .. (351.39,435.63) -- cycle ;
\draw  [draw opacity=0][fill={rgb, 255:red, 155; green, 155; blue, 155 }  ,fill opacity=0.42 ][dash pattern={on 0.84pt off 2.51pt}][line width=0.75]  (410.41,415.39) .. controls (412.84,416.47) and (413.16,421.04) .. (411.13,425.61) .. controls (409.1,430.18) and (405.49,433.01) .. (403.06,431.93) .. controls (400.63,430.85) and (400.31,426.27) .. (402.34,421.7) .. controls (404.36,417.14) and (407.98,414.31) .. (410.41,415.39) -- cycle ;
\draw  [fill={rgb, 255:red, 0; green, 0; blue, 0 }  ,fill opacity=1 ] (407.7,418.96) .. controls (407.73,418.4) and (408.2,417.98) .. (408.75,418) .. controls (409.3,418.03) and (409.73,418.49) .. (409.7,419.05) .. controls (409.68,419.6) and (409.21,420.02) .. (408.66,420) .. controls (408.11,419.98) and (407.68,419.51) .. (407.7,418.96) -- cycle ;
\draw  [fill={rgb, 255:red, 0; green, 0; blue, 0 }  ,fill opacity=1 ] (403.64,428.1) .. controls (403.67,427.54) and (404.14,427.12) .. (404.69,427.14) .. controls (405.24,427.17) and (405.67,427.63) .. (405.64,428.18) .. controls (405.62,428.74) and (405.15,429.16) .. (404.6,429.14) .. controls (404.05,429.11) and (403.62,428.65) .. (403.64,428.1) -- cycle ;
\draw  [fill={rgb, 255:red, 0; green, 0; blue, 0 }  ,fill opacity=1 ] (395.53,446.37) .. controls (395.55,445.82) and (396.02,445.39) .. (396.57,445.42) .. controls (397.12,445.44) and (397.55,445.91) .. (397.52,446.46) .. controls (397.5,447.01) and (397.03,447.44) .. (396.48,447.42) .. controls (395.93,447.39) and (395.5,446.93) .. (395.53,446.37) -- cycle ;
\draw  [fill={rgb, 255:red, 0; green, 0; blue, 0 }  ,fill opacity=1 ] (391.47,455.51) .. controls (391.49,454.96) and (391.96,454.53) .. (392.51,454.56) .. controls (393.06,454.58) and (393.49,455.05) .. (393.46,455.6) .. controls (393.44,456.15) and (392.97,456.58) .. (392.42,456.56) .. controls (391.87,456.53) and (391.44,456.06) .. (391.47,455.51) -- cycle ;
\draw  [draw opacity=0][fill={rgb, 255:red, 155; green, 155; blue, 155 }  ,fill opacity=0.42 ][dash pattern={on 0.84pt off 2.51pt}][line width=0.75]  (398.38,442.46) .. controls (400.81,443.54) and (401.13,448.12) .. (399.11,452.69) .. controls (397.08,457.26) and (393.46,460.09) .. (391.03,459.01) .. controls (388.6,457.93) and (388.28,453.35) .. (390.31,448.78) .. controls (392.34,444.21) and (395.95,441.39) .. (398.38,442.46) -- cycle ;
\draw  [draw opacity=0][fill={rgb, 255:red, 155; green, 155; blue, 155 }  ,fill opacity=0.42 ][dash pattern={on 0.84pt off 2.51pt}][line width=0.75]  (292.74,412.04) .. controls (295.3,411.3) and (298.5,414.59) .. (299.89,419.39) .. controls (301.29,424.19) and (300.35,428.68) .. (297.8,429.42) .. controls (295.24,430.17) and (292.04,426.88) .. (290.65,422.08) .. controls (289.25,417.28) and (290.19,412.79) .. (292.74,412.04) -- cycle ;
\draw [draw opacity=0][fill={rgb, 255:red, 208; green, 2; blue, 27 }  ,fill opacity=0.5 ]   (303.1,419.24) .. controls (323.19,432.79) and (347.45,440.78) .. (342.75,440.53) .. controls (338.05,440.28) and (313.61,441.69) .. (311.48,448.05) .. controls (309.34,454.41) and (316.73,434.44) .. (303.1,419.24) -- cycle ;
\draw  [fill={rgb, 255:red, 0; green, 0; blue, 0 }  ,fill opacity=1 ] (292.98,416.51) .. controls (292.64,416.08) and (292.72,415.45) .. (293.16,415.11) .. controls (293.6,414.78) and (294.23,414.86) .. (294.56,415.29) .. controls (294.9,415.73) and (294.82,416.36) .. (294.38,416.7) .. controls (293.94,417.03) and (293.32,416.95) .. (292.98,416.51) -- cycle ;
\draw  [fill={rgb, 255:red, 0; green, 0; blue, 0 }  ,fill opacity=1 ] (295.77,426.12) .. controls (295.43,425.68) and (295.51,425.05) .. (295.95,424.71) .. controls (296.39,424.38) and (297.02,424.46) .. (297.35,424.9) .. controls (297.69,425.33) and (297.61,425.96) .. (297.17,426.3) .. controls (296.74,426.64) and (296.11,426.55) .. (295.77,426.12) -- cycle ;
\draw  [fill={rgb, 255:red, 0; green, 0; blue, 0 }  ,fill opacity=1 ] (301.35,445.32) .. controls (301.02,444.88) and (301.1,444.26) .. (301.54,443.92) .. controls (301.97,443.58) and (302.6,443.66) .. (302.94,444.1) .. controls (303.28,444.54) and (303.19,445.17) .. (302.76,445.5) .. controls (302.32,445.84) and (301.69,445.76) .. (301.35,445.32) -- cycle ;
\draw  [fill={rgb, 255:red, 0; green, 0; blue, 0 }  ,fill opacity=1 ] (304.15,454.92) .. controls (303.81,454.49) and (303.89,453.86) .. (304.33,453.52) .. controls (304.76,453.18) and (305.39,453.27) .. (305.73,453.7) .. controls (306.07,454.14) and (305.99,454.77) .. (305.55,455.11) .. controls (305.11,455.44) and (304.48,455.36) .. (304.15,454.92) -- cycle ;
\draw [draw opacity=0][fill={rgb, 255:red, 208; green, 2; blue, 27 }  ,fill opacity=0.5 ]   (303.1,419.24) .. controls (322.46,421.67) and (350.56,413.27) .. (345.65,414.69) .. controls (340.74,416.12) and (313.61,441.69) .. (311.48,448.05) .. controls (309.34,454.41) and (316.73,434.44) .. (303.1,419.24) -- cycle ;
\draw  [draw opacity=0][fill={rgb, 255:red, 155; green, 155; blue, 155 }  ,fill opacity=0.42 ][dash pattern={on 0.84pt off 2.51pt}][line width=0.75]  (301.01,440.49) .. controls (303.57,439.75) and (306.77,443.04) .. (308.16,447.84) .. controls (309.56,452.64) and (308.62,457.13) .. (306.07,457.87) .. controls (303.51,458.62) and (300.31,455.33) .. (298.92,450.53) .. controls (297.52,445.73) and (298.46,441.24) .. (301.01,440.49) -- cycle ;
\draw  [dash pattern={on 4.5pt off 4.5pt}]  (284.85,451.91) .. controls (272.85,551.91) and (412.18,563.24) .. (409.18,456.58) ;
\draw  [draw opacity=0][fill={rgb, 255:red, 155; green, 155; blue, 155 }  ,fill opacity=0.42 ][dash pattern={on 0.84pt off 2.51pt}][line width=0.75]  (583.05,432) .. controls (585.71,432) and (587.87,436.05) .. (587.87,441.05) .. controls (587.87,446.05) and (585.71,450.1) .. (583.05,450.1) .. controls (580.4,450.1) and (578.24,446.05) .. (578.24,441.05) .. controls (578.24,436.05) and (580.4,432) .. (583.05,432) -- cycle ;
\draw [draw opacity=0][fill={rgb, 255:red, 208; green, 2; blue, 27 }  ,fill opacity=0.5 ]   (588.96,441.66) .. controls (604.47,460.28) and (625.54,474.73) .. (621.09,473.17) .. controls (616.65,471.61) and (592.79,466.15) .. (588.96,471.66) .. controls (585.14,477.17) and (597.81,460.06) .. (588.96,441.66) -- cycle ;
\draw  [fill={rgb, 255:red, 0; green, 0; blue, 0 }  ,fill opacity=1 ] (582.03,436.36) .. controls (581.83,435.84) and (582.08,435.26) .. (582.6,435.06) .. controls (583.11,434.86) and (583.69,435.11) .. (583.89,435.63) .. controls (584.1,436.14) and (583.84,436.72) .. (583.33,436.92) .. controls (582.81,437.13) and (582.23,436.87) .. (582.03,436.36) -- cycle ;
\draw  [fill={rgb, 255:red, 0; green, 0; blue, 0 }  ,fill opacity=1 ] (582.03,446.36) .. controls (581.83,445.84) and (582.08,445.26) .. (582.6,445.06) .. controls (583.11,444.86) and (583.69,445.11) .. (583.89,445.63) .. controls (584.1,446.14) and (583.84,446.72) .. (583.33,446.92) .. controls (582.81,447.13) and (582.23,446.87) .. (582.03,446.36) -- cycle ;
\draw  [fill={rgb, 255:red, 0; green, 0; blue, 0 }  ,fill opacity=1 ] (582.03,466.36) .. controls (581.83,465.84) and (582.08,465.26) .. (582.6,465.06) .. controls (583.11,464.86) and (583.69,465.11) .. (583.89,465.63) .. controls (584.1,466.14) and (583.84,466.72) .. (583.33,466.92) .. controls (582.81,467.13) and (582.23,466.87) .. (582.03,466.36) -- cycle ;
\draw  [fill={rgb, 255:red, 0; green, 0; blue, 0 }  ,fill opacity=1 ] (582.03,476.36) .. controls (581.83,475.84) and (582.08,475.26) .. (582.6,475.06) .. controls (583.11,474.86) and (583.69,475.11) .. (583.89,475.63) .. controls (584.1,476.14) and (583.84,476.72) .. (583.33,476.92) .. controls (582.81,477.13) and (582.23,476.87) .. (582.03,476.36) -- cycle ;
\draw [draw opacity=0][fill={rgb, 255:red, 208; green, 2; blue, 27 }  ,fill opacity=0.5 ]   (588.96,441.66) .. controls (606.87,449.39) and (636.2,449.17) .. (631.09,449.17) .. controls (625.98,449.17) and (592.79,466.15) .. (588.96,471.66) .. controls (585.14,477.17) and (597.81,460.06) .. (588.96,441.66) -- cycle ;
\draw  [draw opacity=0][fill={rgb, 255:red, 155; green, 155; blue, 155 }  ,fill opacity=0.42 ][dash pattern={on 0.84pt off 2.51pt}][line width=0.75]  (583.05,461.63) .. controls (585.71,461.63) and (587.87,465.68) .. (587.87,470.68) .. controls (587.87,475.68) and (585.71,479.73) .. (583.05,479.73) .. controls (580.4,479.73) and (578.24,475.68) .. (578.24,470.68) .. controls (578.24,465.68) and (580.4,461.63) .. (583.05,461.63) -- cycle ;
\draw  [draw opacity=0][fill={rgb, 255:red, 155; green, 155; blue, 155 }  ,fill opacity=0.42 ][dash pattern={on 0.84pt off 2.51pt}][line width=0.75]  (642.07,441.39) .. controls (644.5,442.47) and (644.83,447.04) .. (642.8,451.61) .. controls (640.77,456.18) and (637.16,459.01) .. (634.73,457.93) .. controls (632.3,456.85) and (631.97,452.27) .. (634,447.7) .. controls (636.03,443.14) and (639.64,440.31) .. (642.07,441.39) -- cycle ;
\draw  [fill={rgb, 255:red, 0; green, 0; blue, 0 }  ,fill opacity=1 ] (639.37,444.96) .. controls (639.4,444.4) and (639.86,443.98) .. (640.41,444) .. controls (640.97,444.03) and (641.39,444.49) .. (641.37,445.05) .. controls (641.34,445.6) and (640.88,446.02) .. (640.33,446) .. controls (639.77,445.98) and (639.35,445.51) .. (639.37,444.96) -- cycle ;
\draw  [fill={rgb, 255:red, 0; green, 0; blue, 0 }  ,fill opacity=1 ] (635.31,454.1) .. controls (635.34,453.54) and (635.8,453.12) .. (636.35,453.14) .. controls (636.91,453.17) and (637.33,453.63) .. (637.31,454.18) .. controls (637.29,454.74) and (636.82,455.16) .. (636.27,455.14) .. controls (635.71,455.11) and (635.29,454.65) .. (635.31,454.1) -- cycle ;
\draw  [fill={rgb, 255:red, 0; green, 0; blue, 0 }  ,fill opacity=1 ] (627.19,472.37) .. controls (627.22,471.82) and (627.68,471.39) .. (628.24,471.42) .. controls (628.79,471.44) and (629.22,471.91) .. (629.19,472.46) .. controls (629.17,473.01) and (628.7,473.44) .. (628.15,473.42) .. controls (627.6,473.39) and (627.17,472.93) .. (627.19,472.37) -- cycle ;
\draw  [fill={rgb, 255:red, 0; green, 0; blue, 0 }  ,fill opacity=1 ] (623.13,481.51) .. controls (623.16,480.96) and (623.62,480.53) .. (624.18,480.56) .. controls (624.73,480.58) and (625.16,481.05) .. (625.13,481.6) .. controls (625.11,482.15) and (624.64,482.58) .. (624.09,482.56) .. controls (623.54,482.53) and (623.11,482.06) .. (623.13,481.51) -- cycle ;
\draw  [draw opacity=0][fill={rgb, 255:red, 155; green, 155; blue, 155 }  ,fill opacity=0.42 ][dash pattern={on 0.84pt off 2.51pt}][line width=0.75]  (630.05,468.46) .. controls (632.48,469.54) and (632.8,474.12) .. (630.77,478.69) .. controls (628.74,483.26) and (625.13,486.09) .. (622.7,485.01) .. controls (620.27,483.93) and (619.95,479.35) .. (621.97,474.78) .. controls (624,470.21) and (627.62,467.39) .. (630.05,468.46) -- cycle ;
\draw  [draw opacity=0][fill={rgb, 255:red, 155; green, 155; blue, 155 }  ,fill opacity=0.42 ][dash pattern={on 0.84pt off 2.51pt}][line width=0.75]  (524.41,438.04) .. controls (526.96,437.3) and (530.16,440.59) .. (531.56,445.39) .. controls (532.95,450.19) and (532.02,454.68) .. (529.46,455.42) .. controls (526.91,456.17) and (523.71,452.88) .. (522.31,448.08) .. controls (520.92,443.28) and (521.86,438.79) .. (524.41,438.04) -- cycle ;
\draw [draw opacity=0][fill={rgb, 255:red, 208; green, 2; blue, 27 }  ,fill opacity=0.5 ]   (534.77,445.24) .. controls (554.86,458.79) and (579.12,466.78) .. (574.42,466.53) .. controls (569.71,466.28) and (545.27,467.69) .. (543.14,474.05) .. controls (541.01,480.41) and (548.4,460.44) .. (534.77,445.24) -- cycle ;
\draw  [fill={rgb, 255:red, 0; green, 0; blue, 0 }  ,fill opacity=1 ] (524.65,442.51) .. controls (524.31,442.08) and (524.39,441.45) .. (524.83,441.11) .. controls (525.26,440.78) and (525.89,440.86) .. (526.23,441.29) .. controls (526.57,441.73) and (526.49,442.36) .. (526.05,442.7) .. controls (525.61,443.03) and (524.98,442.95) .. (524.65,442.51) -- cycle ;
\draw  [fill={rgb, 255:red, 0; green, 0; blue, 0 }  ,fill opacity=1 ] (527.44,452.12) .. controls (527.1,451.68) and (527.18,451.05) .. (527.62,450.71) .. controls (528.06,450.38) and (528.68,450.46) .. (529.02,450.9) .. controls (529.36,451.33) and (529.28,451.96) .. (528.84,452.3) .. controls (528.4,452.64) and (527.77,452.55) .. (527.44,452.12) -- cycle ;
\draw  [fill={rgb, 255:red, 0; green, 0; blue, 0 }  ,fill opacity=1 ] (533.02,471.32) .. controls (532.68,470.88) and (532.76,470.26) .. (533.2,469.92) .. controls (533.64,469.58) and (534.27,469.66) .. (534.6,470.1) .. controls (534.94,470.54) and (534.86,471.17) .. (534.42,471.5) .. controls (533.99,471.84) and (533.36,471.76) .. (533.02,471.32) -- cycle ;
\draw  [fill={rgb, 255:red, 0; green, 0; blue, 0 }  ,fill opacity=1 ] (535.81,480.92) .. controls (535.47,480.49) and (535.56,479.86) .. (535.99,479.52) .. controls (536.43,479.18) and (537.06,479.27) .. (537.4,479.7) .. controls (537.73,480.14) and (537.65,480.77) .. (537.21,481.11) .. controls (536.78,481.44) and (536.15,481.36) .. (535.81,480.92) -- cycle ;
\draw [draw opacity=0][fill={rgb, 255:red, 208; green, 2; blue, 27 }  ,fill opacity=0.5 ]   (534.77,445.24) .. controls (554.12,447.67) and (582.23,439.27) .. (577.32,440.69) .. controls (572.41,442.12) and (545.27,467.69) .. (543.14,474.05) .. controls (541.01,480.41) and (548.4,460.44) .. (534.77,445.24) -- cycle ;
\draw  [draw opacity=0][fill={rgb, 255:red, 155; green, 155; blue, 155 }  ,fill opacity=0.42 ][dash pattern={on 0.84pt off 2.51pt}][line width=0.75]  (532.68,466.49) .. controls (535.23,465.75) and (538.43,469.04) .. (539.83,473.84) .. controls (541.22,478.64) and (540.29,483.13) .. (537.73,483.87) .. controls (535.18,484.62) and (531.98,481.33) .. (530.59,476.53) .. controls (529.19,471.73) and (530.13,467.24) .. (532.68,466.49) -- cycle ;
\draw  [dash pattern={on 4.5pt off 4.5pt}]  (516.52,477.91) .. controls (504.52,577.91) and (643.85,589.24) .. (640.85,482.58) ;
\draw  [draw opacity=0][fill={rgb, 255:red, 126; green, 211; blue, 33 }  ,fill opacity=0.7 ] (505.25,397.72) .. controls (505.25,391.65) and (514.88,386.72) .. (526.75,386.72) .. controls (538.63,386.72) and (548.25,391.65) .. (548.25,397.72) .. controls (548.25,403.8) and (538.63,408.72) .. (526.75,408.72) .. controls (514.88,408.72) and (505.25,403.8) .. (505.25,397.72) -- cycle ;
\draw  [draw opacity=0][fill={rgb, 255:red, 126; green, 211; blue, 33 }  ,fill opacity=0.7 ] (561.2,380.92) .. controls (561.2,374.84) and (570.82,369.92) .. (582.7,369.92) .. controls (594.57,369.92) and (604.2,374.84) .. (604.2,380.92) .. controls (604.2,386.99) and (594.57,391.92) .. (582.7,391.92) .. controls (570.82,391.92) and (561.2,386.99) .. (561.2,380.92) -- cycle ;
\draw  [draw opacity=0][fill={rgb, 255:red, 208; green, 2; blue, 27 }  ,fill opacity=0.5 ] (524.26,411.42) -- (528.82,435.33) -- (520.48,435.33) -- cycle ;
\draw  [draw opacity=0][fill={rgb, 255:red, 208; green, 2; blue, 27 }  ,fill opacity=0.5 ] (583.7,391.92) -- (587.54,424.42) -- (580.51,424.42) -- cycle ;
\draw  [draw opacity=0][fill={rgb, 255:red, 208; green, 2; blue, 27 }  ,fill opacity=0.5 ] (638.26,410.92) -- (642.82,434.83) -- (634.48,434.83) -- cycle ;
\draw  [draw opacity=0][fill={rgb, 255:red, 126; green, 211; blue, 33 }  ,fill opacity=0.7 ] (612.36,398.22) .. controls (612.36,392.15) and (621.99,387.22) .. (633.86,387.22) .. controls (645.74,387.22) and (655.36,392.15) .. (655.36,398.22) .. controls (655.36,404.3) and (645.74,409.22) .. (633.86,409.22) .. controls (621.99,409.22) and (612.36,404.3) .. (612.36,398.22) -- cycle ;
\draw   (232.33,476.45) -- (248.37,476.45) -- (248.37,472.97) -- (259.06,479.94) -- (248.37,486.91) -- (248.37,483.42) -- (232.33,483.42) -- cycle ;
\draw   (448.33,476.45) -- (464.37,476.45) -- (464.37,472.97) -- (475.06,479.94) -- (464.37,486.91) -- (464.37,483.42) -- (448.33,483.42) -- cycle ;
\draw  [draw opacity=0][fill={rgb, 255:red, 208; green, 2; blue, 27 }  ,fill opacity=0.5 ] (503.24,456.67) -- (525.27,467.02) -- (520.38,473.77) -- cycle ;
\draw  [draw opacity=0][fill={rgb, 255:red, 126; green, 211; blue, 33 }  ,fill opacity=0.7 ] (467.75,444.42) .. controls (467.75,438.34) and (477.38,433.42) .. (489.25,433.42) .. controls (501.13,433.42) and (510.75,438.34) .. (510.75,444.42) .. controls (510.75,450.49) and (501.13,455.42) .. (489.25,455.42) .. controls (477.38,455.42) and (467.75,450.49) .. (467.75,444.42) -- cycle ;
\draw  [draw opacity=0][fill={rgb, 255:red, 208; green, 2; blue, 27 }  ,fill opacity=0.5 ] (656.6,456.03) -- (639.92,473.77) -- (635.03,467.02) -- cycle ;
\draw  [draw opacity=0][fill={rgb, 255:red, 126; green, 211; blue, 33 }  ,fill opacity=0.7 ] (653.86,444.4) .. controls (653.86,438.32) and (663.49,433.4) .. (675.36,433.4) .. controls (687.24,433.4) and (696.86,438.32) .. (696.86,444.4) .. controls (696.86,450.48) and (687.24,455.4) .. (675.36,455.4) .. controls (663.49,455.4) and (653.86,450.48) .. (653.86,444.4) -- cycle ;
\draw  [draw opacity=0][fill={rgb, 255:red, 208; green, 2; blue, 27 }  ,fill opacity=0.5 ] (583.74,501.73) -- (579.48,481.83) -- (587.24,481.83) -- cycle ;
\draw  [draw opacity=0][fill={rgb, 255:red, 126; green, 211; blue, 33 }  ,fill opacity=0.7 ] (563.2,515.42) .. controls (563.2,509.34) and (572.6,504.42) .. (584.2,504.42) .. controls (595.8,504.42) and (605.2,509.34) .. (605.2,515.42) .. controls (605.2,521.49) and (595.8,526.42) .. (584.2,526.42) .. controls (572.6,526.42) and (563.2,521.49) .. (563.2,515.42) -- cycle ;


\draw (477,434) node [anchor=north west][inner sep=0.75pt]  [font=\scriptsize] [align=left] {$Z_{\ell}^{(3)}$};
\draw (662,434) node [anchor=north west][inner sep=0.75pt]  [font=\scriptsize] [align=left] {$Z_{\ell}^{(3)}$};
\draw (570.55,504.5) node [anchor=north west][inner sep=0.75pt]  [font=\scriptsize] [align=left] {$Z_{\ell}^{(3)}$};

\draw (117.36,583) node [anchor=north west][inner sep=0.75pt]  [font=\scriptsize] [align=left] {$G_{\ell}^1\cong Z_{\ell}^{(3)}$};
\draw (335.19,583) node [anchor=north west][inner sep=0.75pt]  [font=\scriptsize] [align=left] {$G_{\ell}^1(2)$};
\draw (513,389) node [anchor=north west][inner sep=0.75pt]  [font=\scriptsize] [align=left] {$Z_{\ell}^{(3)}$};
\draw (570,370) node [anchor=north west][inner sep=0.75pt]  [font=\scriptsize] [align=left] {$Z_{\ell}^{(3)}$};
\draw (620,389) node [anchor=north west][inner sep=0.75pt]  [font=\scriptsize] [align=left] {$Z_{\ell}^{(3)}$};
\draw (573,583) node [anchor=north west][inner sep=0.75pt]  [font=\scriptsize] [align=left] {$G_{\ell}^2$};

\end{tikzpicture}
\caption{Construction of $G_{\ell}^2$ from $G_{\ell}^1$ for $k=3$}\label{FIGHRE:G_L^R}
\end{figure}

    In the following we argue that these~$G_{\ell}^r$ satisfy property~\eqref{it:upperbound}.
    For this let~$\varepsilon>0$ be given.
    Clearly, we may assume that~$\varepsilon\ll 1$.
    By property~\eqref{it:upperbound} (and induction) and Lemma~\ref{lem:zyclesupper}, there are~$\ell_1,
    \ell_2\in\mathds{N}$ such that~$\gamma\big(G_{m\cdot\ell_1}^{r-1}\big)\leq\frac{r-1}{r}$ holds for every~$m\in\mathds{N}$ and~$\gamma\big(Z_{\ell_2}^{(k)}\big)\leq\varepsilon/2$.
    Setting~$\ell=\ell_1\cdot\ell_2$, we have~$\gamma\big(G_{m\cdot\ell}^{r-1}\big)\leq\frac{r-1}{r}$ and, using the same blow-up argument as above,~$\gamma\big(Z_{m\cdot\ell}^{(k)}\big)\leq\gamma\big(Z_{\ell_2}^{(k)}\big)\leq\varepsilon/2$ for every~$m\in\mathds{N}$.
    
    Given~$m\in\mathds{N}$, let~$n\in\mathds{N}$ be sufficiently large and let~$H$ be a~$k$-graph on~$n$ vertices with~$\delta(H)\geq(\frac{r-1}{r}+\varepsilon)n$.
    Due to Proposition~\ref{prop:supersaturation},~$\gamma\left(G_{m\cdot\ell}^{r-1}(k-1)\right)=\gamma\left(G_{m\cdot\ell}^{r-1}\right)$, whence our choice of~$\ell$ yields~$G_{m\cdot\ell}^{r-1}(k-1)\subseteq H$.
    Note that the minimum codegree of~$H$ implies that for any~$(k-1)$-sets~$p_1,\dots,p_{r-1}\in V(H)^{(k-1)}$, we have $$\Big\vert\bigcap_{i\in[r-1]}N(p_i)\Big\vert\geq\Big(\frac{1}{r}+\varepsilon\Big) n\,.$$
    Using the minimum codegree of~$H$ once more, this entails that~$\delta\big(H\big[\bigcap_{i\in[r-1]}N(p_i)\big]\big)\geq\varepsilon n$.
    In fact, for every set~$X\subseteq V(H)$ with~$\vert X\vert\leq \varepsilon n/4$, we have $$\delta\Big(H\big[\bigcap_{i\in[r-1]}N(p_i)\setminus X\big]\Big)\geq\frac{3}{4}\varepsilon n\,.$$
    Therefore, the choices of~$\ell$ and~$n$ imply that~$Z_{m\cdot\ell}^{(k)}\subseteq H\big[\bigcap_{i\in[r-1]}N(p_i)\setminus X\big]$.
    In particular, this allows us to find copies of~$Z_{m\cdot\ell}^{(k)}$ in each~$\bigcap_{v\in S}N(x_v^1\cdots x_v^{k-1})$ with~$S\in V(G_{m\cdot\ell}^{r-1})^{(r-1)}$ that avoid all previously used vertices since at any point at most~$\varepsilon n/4$ vertices have been used\footnote{Of course, due to Proposition~\ref{prop:supersaturation}, we really only need to find a homomorphism from~$G_{m\cdot\ell}^r$ into~$H$.}.
    Thus a copy of~$G_{m\cdot\ell}^r$ is obtained in~$H$.

    Next, we prove that the~$k$-graphs~$G_{\ell}^r$ satisfy property~\eqref{it:lowerbound}.
    Choose~$p_0\gg\ell,r,k$ (in particular,~$p_0$ is at least as large as the~$p_0$ that is guaranteed when appealing to property~\eqref{it:lowerbound} with~$r-1$ and~$\ell$) and let~$p\geq p_0$ be a prime,~$\eta>0$, and~$n\in\mathds{N}$ with~$p\mid n-(r-1)\cdot\big\lfloor\frac{(1-\eta)n}{r}\big\rfloor$.
    Assume for the sake of contradiction that~$G^r_{\ell}\subseteq H^r_{\mathds{F}_p}(\eta,n)$.
    First, observe that for every~$i\in[r-1]$, there must be some~$v\in V(G_{\ell}^{r-1})$ such that 
    all the vertices $x_v^1,x_v^2,\ldots,x_v^{k-1}$ are contained in~$V_i$.
    Indeed, if this was not the case for some~$i\in[r-1]$, then for every~$v\in V(G_{\ell}^{r-1})$, one of the copies of~$v$ is contained in~$\bigcup_{j\in[r]\setminus\{i\}}V_j$.
    In other words, there is a copy of~$G_{\ell}^{r-1}$ contained in~$H^r_{\mathds{F}_p}(\eta,n)-V_i$.
    Since~$H^r_{\mathds{F}_p}(\eta,n)-V_i$ is isomorphic to~$H^{r-1}_{\mathds{F}_p}(\eta',n')$ with $$n'=n-\Big\lfloor\frac{(1-\eta)n}{r}\Big\rfloor\text{ and }\eta'=1-\frac{(r-1)\big\lfloor\frac{(1-\eta)n}{r}\big\rfloor}{n'}$$ (note that for this choice we in particular have~$p\mid n'-(r-2)\cdot\big\lfloor\frac{(1-\eta')n'}{r-1}\big\rfloor$), this contradicts that by property~\eqref{it:lowerbound} (and induction),~$G_{\ell}^{r-1}\not\subseteq H^{r-1}_{\mathds{F}_p}(\eta',n')$.

    Thus, we may assume that there are~$v_1,\dots,v_{r-1}\in V(G_{\ell}^{r-1})$ such that for every~$i\in[r-1]$, we have~$x_{v_i}^1,\ldots,x_{v_i}^{k-1}\in V_i$.
    Setting~$S=\{v_1,\dots,v_{r-1}\}$, the induced~$k$-graph $$H^r_{\mathds{F}_p}(\eta,n)\Big[\bigcap_{i\in[r-1]}N(x_{v_i}^1\cdots x_{v_i}^{k-1})\Big]$$ must contain a copy of~$Z_{\ell}^{(k)}$.
    But note that $\bigcap_{i\in[r-1]}N(x_{v_i}^1\cdots x_{v_i}^{k-1})=V_r$ and that $H^r_{\mathds{F}_p}(\eta,n)[V_r]\cong\mathds{F}^{(k)}_p(|V_r|)$.
    This entails a contradiction since by Lemma~\ref{lem:Fpprop} and the choice of~$p$,~$\mathds{F}^{(k)}_p(|V_r|)$ does not contain a copy of~$Z_{\ell}^{(k)}$.

\section*{Acknowledgements}
The third author is very grateful for the outstanding hospitality of Guanghui Wang and his group at Shandong University.
We thank David Conlon, Sim\'on Piga, and Haotian Yang for helpful discussions.
Research was supported by the Postdoctoral Fellowship Program of CPSF under Grant Number GZC20252020 and the China Postdoctoral Science Foundation 2025M783118 (Weichan Liu), the Young Scientist Fellowship IBS-R029-Y7 (Bjarne Sch\"ulke), and the China Postdoctoral Science Foundation 2023M742092 (Wanting Sun), and Shandong University.

\end{document}